\documentclass[11pt]{article}
\usepackage[utf8]{inputenc}
\usepackage[T1]{fontenc}
\usepackage{amsmath,amssymb,graphicx}
\usepackage[a4paper,margin=30mm]{geometry}
\usepackage[pdfborder={0 0 0}]{hyperref}
\title{On the Number of Neighborly Simplices in \(\mathbb{R}^d\)}
\author{Andrzej P. Kisielewicz}
\date{}
\usepackage{pgf}
\definecolor{figcoloraa}{rgb}{0.65882,0.65882,0.66275}
\definecolor{figcolorab}{rgb}{0.07451,0.078431,0.07451}
\definecolor{figcolorac}{rgb}{0.54902,0.54902,0.54902}
\definecolor{figcolorad}{rgb}{0.91373,0.91765,0.91765}
\definecolor{figcolorae}{rgb}{0.76863,0.77255,0.77255}
\newcommand{\InflationFigure}{%
}

\begin{document}
\maketitle
\begin{abstract}

Two \textit{d}-dimensional simplices in \(\mathbb {R}^d\) are neighborly if its intersection is a \((d-1)\)-dimensional set. A family of \textit{d}-dimensional simplices in \(\mathbb {R}^d\) is called neighborly if every two simplices of the family are neighborly. Let \(S_d\) be the maximal cardinality of a neighborly family of \textit{d}-dimensional simplices in \(\mathbb {R}^d\). Based on the structure of some codes \(V\subset \{0,1,*\}^n\) it is shown that \(\lim _{d\rightarrow \infty }(2^{d+1}-S_d)=\infty \). Moreover, a result on the structure of codes \(V\subset \{0,1,*\}^n\) is given.

\end{abstract}


\hypertarget{Sec1}{}\section{Introduction}

Let \(A=\{0,1,*\}\), and let \(A^n\) be the set of all words \(v=v_1 \ldots v_n\) over the alphabet \textit{A}, that is, \(v_i\in A\) for \(i\in [n]=\{1, \ldots ,n\}\). Two words \(v,u\in A^n\) are called \textit{dichotomous} (\textit{at the i-th position}) if \(v_i+u_i=1\) (\(u_i,v_i\in \{0,1\}\)) for some \(i\in [n]\) (compare [\hyperlink{ref-CR11}{11}, Section 10]) and they are called \textit{neighborly} if there is precisely one such \textit{i}. Two neighborly words \(v,w\in A^n\) are a \textit{twin pair} if \(w_j=v_j\) for all \(j\in [n]{\setminus } \{i\}\). A family of words \(V\subset A^n\) is called a \textit{( dichotomous) code} if every two words in \textit{V} are dichotomous. A code \(V\subset A^n\) is a \textit{d-code} if \(|{\text {prop}}\, (v)|=|\{i\in [n]:v_i\ne *\}|=d\ge 1\) for every \(v\in V\). A family of words \(V\subset A^n\) is called a \textit{neighborly} code if every two words in \textit{V} are neighborly. By \(M_d\) we denote the maximal cardinality of a neighborly \textit{d}-code without twin pairs.

Two \textit{d}-dimensional simplices in \(\mathbb {R}^d\) are \textit{neighborly} if its intersection is a \((d-1)\)-dimensional set. A family of \textit{d}-dimensional simplices in \(\mathbb {R}^d\) is called \textit{neighborly} if every two simplices of the family are neighborly. Let \(S_d\) be the maximal cardinality of a neighborly family of \textit{d}-dimensional simplices in \(\mathbb {R}^d\). A long standing conjecture says that \(S_d=2^d\) [\hyperlink{ref-CR5}{5}]. It is verified for dimensions \(d\le 3\). In [\hyperlink{ref-CR14}{14}] J. Zaks showed that \(S_3=8\) (earlier V. Baston proved in [\hyperlink{ref-CR6}{6}] that \(S_3\le 9\)), and in [\hyperlink{ref-CR13}{13}] that \(S_d\ge 2^d\). M. Perles proved the estimation \(S_d\le 2^{d+1}\) [\hyperlink{ref-CR12}{12}], and M. Aigner and G. Ziegler showed that \(S_d\le 2^{d+1}-1\) [\hyperlink{ref-CR1}{1}, Chapter 14]. Recently in [\hyperlink{ref-CR9}{9}, \hyperlink{ref-CR10}{10}] it was shown that \(S_d\le 2^{d+1}-2\). In this note we prove

\hypertarget{FPar1}{}\paragraph{Theorem 1.}

If \(S_d\) is the maximal cardinality of a neighborly family of \textit{d}-dimensional simplices in \(\mathbb {R}^d\), then

\hypertarget{Equ13}{}\begin{equation*}
\begin{aligned} \lim _{d\rightarrow \infty }(2^{d+1}-S_d)=\infty . \end{aligned}
\end{equation*}

Below we describe a passing from neighborly simplices to neighborly \textit{d}-codes.

Let \textit{S} be a neighborly family of \textit{d}-dimensional simplices in \(\mathbb {R}^d\), and let \(H_1, \ldots ,H_n\) denote all different hyperplanes spanned by facets of simplices in \textit{S} (compare Fig. \hyperlink{Fig2}{2} in Sect. \hyperlink{Sec3}{3}). Let \(H^0_i,H^1_i\) be two sides of \(H_i\), \(i \in [n]\). Fix \(\sigma _d\in S\) and define a word \(v=v(\sigma _d) \in A^n\):

\hypertarget{Equ14}{}\begin{equation*}
\begin{aligned} v_i= {\left\{ \begin{array}{ll} 0, &{} \text {if }H_i\text { is spanned by a facet of }\sigma _d\text { and }\sigma _d\subset H_i^0,\\ 1, &{} \text {if }H_i\text { is spanned by a facet of }\sigma _d\text { and }\sigma _d\subset H_i^1,\\ *, &{} \text {otherwise.} \end{array}\right. } \end{aligned}
\end{equation*}

Let \(V=\{v(\sigma _d):\sigma _d \in S\}\). In the same way as in [\hyperlink{ref-CR6}{6}, Lemmas 1-4] (for the case \(d=3\)) we show that \textit{V} is a neighborly \((d+1)\)-code without twin pairs (compare also [\hyperlink{ref-CR1}{1}, Chapter 14]). Since \(|S|=|V|\), we have \(S_d\le M_{d+1}\). Therefore, to prove Theorem \hyperlink{FPar1}{1} we shall prove

\hypertarget{FPar2}{}\paragraph{Theorem 2.}

If \(M_d\) is the maximal cardinality of a neighborly \textit{d}-code without twin pairs, then

\hypertarget{Equ15}{}\begin{equation*}
\begin{aligned} \lim _{d\rightarrow \infty }(2^{d}-M_d)=\infty . \end{aligned}
\end{equation*}

A geometric interpretation of a word \(v=v_1\ldots v_n\in A^n\) is the box \(\breve{v}=\breve{v}_1\times \cdots \times \breve{v}_n\subset [0,2]^n\), called a \textit{realization} of \textit{v}, where \(\breve{v}_i=[0,1)\) if \(v_i=0\), \(\breve{v}_i=[1,2]\) if \(v_i=1\) and \(\breve{v}_i=[0,2]\) if \(v_i=*\) (compare Fig. \hyperlink{Fig1}{1}). If \(V\subset A^n\) is a code, then the boxes of the family \(\breve{V}=\{\breve{v}:v\in V\}\) are mutually disjoint. Because of this interpretation we shall use the following notation: \(|v|=2^{n-p}\), where \(p=|{\text {prop}}\, (v)|\) and \({\text {vol}}\, (V)=\sum _{v\in V} |v|\) if \(V\subset A^n\) is a code. (Clearly, \({\text {vol}}\, (V)=m_n(\bigcup _{v\in V}\breve{v})\), where \(m_n\) is the \textit{n}-dimensional Lebesgue measure.) To simplify notations, throughout the paper we shall working mainly with words \(v\in A^n\) rather than boxes \(\breve{v}\subset [0,2]^n\). However, it is very useful to keep in mind the above geometric interpretation of codes \(V\subset A^n\) as it makes reasoning easier. For example, by this interpretation it is immediately clear that \(|V|\le 2^d\) for every \textit{d}-code \(V\subset A^n\) (compare [\hyperlink{ref-CR4}{4}, Theorem 1] for a more general case). Indeed, since \textit{V} is a \textit{d}-code, we have \(m_n(\breve{v})=2^{n-d}\) for every \(v\in V\). As \(\breve{v}\subset [0,2]^n\) for \(\breve{v}\in \breve{V}\) and boxes in the family \(\breve{V}=\{\breve{v}:v\in V\}\) are mutually disjoint, it follows that \(m_n(\breve{V})=2^{n-d}|V|\le 2^n\) and thus, \(|V|\le 2^d\). Hence, \(M_d\le 2^d\). A slightly better upper bound of \(M_d\) is given in [\hyperlink{ref-CR9}{9}, \hyperlink{ref-CR10}{10}, Theorem 1] were we showed that

\hypertarget{Equ1}{}\begin{equation*}
\begin{aligned} M_d\le 2^d-2. \end{aligned}\tag{1}
\end{equation*}

Our proof of Theorem \hyperlink{FPar2}{2} is based on properties of neighborly \textit{d}-code without twin pairs. This technique was introduced by V. Baston in [\hyperlink{ref-CR6}{6}], and next it was used by J. Zaks and M. Perles in [\hyperlink{ref-CR12}{12}, \hyperlink{ref-CR14}{14}]. Originally Baston considered families of strings from the set \(\{-1,0,1\}^n\) arranged as rows of a \textit{matrix representation} of neighborly family of simplices (the translation into our notation is as follows: \(-1=1,1=0\) and \(0=*\) and rows of a matrix representation form neighborly \textit{d}-code without twin pairs). He used combinatorial properties of such matrices and its relationships with neighborly simplices. In [\hyperlink{ref-CR14}{14}], Zaks used the machinery introduced by Baston together with tools from graph theory (Graham--Pollak theorem) as well as a computer support. Our approach is heavily related to a geometrical interpretation of a neighborly \textit{d}-code as a set of boxes and we do not use, unlike Baston and Zaks, any relationships between neighborly codes and neighborly simplices that generated such codes.

A neighborly \textit{d}-code without twin pairs is a very special case of a more general set of words \(V\subset A^n\) which is called a \textit{k-neighborly family} in which every two words from \textit{V} differ by 0 and 1 in at least one and at most \(k\in [n]\) positions [\hyperlink{ref-CR2}{2}, \hyperlink{ref-CR3}{3}, \hyperlink{ref-CR7}{7}]. Neighborly families are closely related to Graham--Pollak theorem, while \textit{k}-neighborly families are related to coverings of complete graphs by bicliques [\hyperlink{ref-CR2}{2}].

\hypertarget{Sec2}{}\section{The Structure of Neighborly Codes}

In this section, we give two results on the structure of neighborly \textit{d}-codes.

Let \(V\subset A^n\) be a code, \(j\in [n]\), and let

\hypertarget{Equ16}{}\begin{equation*}
\begin{aligned} V^{j,a}=\{v\in V:v_j=a\}, \end{aligned}
\end{equation*}

where \(a \in \{0,1,*\}\). If \textit{V} is neighborly, then every two words \(v,u\in V\) are dichotomous at precisely one position \(i\in [n]\). This property enforces a certain structure of \textit{V} which is described in the following lemma:

\hypertarget{FPar3}{}\paragraph{Lemma 1.}

Let \(V\subset A^n\) be a neighborly code, \(j\in [n]\), and let

\hypertarget{Equ17}{}\begin{equation*}
\begin{aligned} C^j_0=\{k\in [n]{\setminus }\{j\}:\mathrm{there\,\, are}\,\, v,u\in V^{j,0} \,\, \mathrm{such \,\, that}\,\, v_k+u_k=1\} \end{aligned}
\end{equation*}

and

\hypertarget{Equ18}{}\begin{equation*}
\begin{aligned} C^j_1=\{k\in [n]{\setminus }\{j\}:\mathrm{there\,\, are}\,\, v,u\in V^{j,1} \,\, \mathrm{such \,\, that}\,\, v_k+u_k=1\}. \end{aligned}
\end{equation*}

If \(|V^{j,0}|\ge 2\), then \(C^j_0\ne \emptyset \) and if \(|V^{j,1}|\ge 2\), then \( C^j_1\ne \emptyset \). Moreover, if \(C^j_0\ne \emptyset \), then for every \(v\in V^{j,1}\) and every \(k\in C^j_0\) we have \(v_k=*\). Similarly, if \(C^j_1\ne \emptyset \), then for every \(v\in V^{j,0}\) and every \(k\in C^j_1\) we have \(v_k=*\). In particular, \(C^j_0\cap C^j_1=\emptyset \).

\hypertarget{FPar4}{}\paragraph{Proof.}

Let \(v,u\in V^{j,0}\) be two distinct words. Since \(v_j=u_j=0\) and \textit{V} is a code, there is \(k\in [n]{\setminus } \{j\}\) such that \(v_k+u_k=1\). Thus, \(k\in C^j_0\). In the same way we show that \(C^j_1\ne \emptyset \) if \(|V^{j,1}|\ge 2\).

Let \(v\in V^{j,1}\) and \(k\in C^j_0\). Since \(v_j=1\), the word \textit{v} is dichotomous with every word in \(V^{j,0}\) at the \textit{j}-th position. Let \(u,w\in V^{j,0}\) be such that \(u_k+w_k=1\). Then \(u_k=0, w_k=1\) or \(u_k=1, w_k=0\). Thus, if \(v_k\in \{0,1\}\), then the words \textit{v}, \textit{u} or \textit{v}, \textit{w} are dichotomous at the \textit{j}-th and the \textit{k}-th position which is impossible. Hence \(v_k=*\). In the same way we consider the case \(v\in V^{j,0}\) and \(k\in C^j_1\).

It follows from the above that \(C^j_0\cap C^j_1=\emptyset \), as if \(k\in C^j_0\cap C^j_1\), then there are four words \textit{v}, \textit{u}, \textit{p}, \textit{q} such that \(v,u\in V^{j,0}\) and \(p,q\in V^{j,1}\) where \(v_k,u_k,p_k,q_k\in \{0,1\}\) which is, as we showed above, impossible. \(\square \)

Let \(v \in A^n\), \(\sigma \) be a permutation of the set [\textit{n}] and let \(\bar{\sigma }(v)=v_{\sigma (1)}\ldots v_{\sigma (n)}\). For \(i\in [n]\) let \(h_i: A\rightarrow A\) be such that \(h_i=f\), where \(f(0)=1, f(1)=0\) and \(f(*)=*\) or \(h_i=\textrm{id}\). The function \textit{f} is called a \textit{flip}. Let \(h: A^n\rightarrow A^n\) be defined be the formula \(h(v)=h_1(v_1)\ldots h_n(v_n)\). Two codes \(V,U \subset A^n\) are \textit{isomorphic} if there is a function \(h\circ \bar{\sigma }\) such that \(U=h\circ \bar{\sigma }(V)\).

Similarly as in [\hyperlink{ref-CR6}{6}], words of a code \(V\subset A^n\) can be represented as rows of a \(|V|\times n\) matrix \textit{M}(\textit{V}). Thus, two codes \(V,U\subset A^n\) are isomorphic if there are a permutation of columns and rows in \textit{M}(\textit{V}) and flips \textit{f}(\textit{a}), \(a\in A\), of letters in some columns of \textit{M}(\textit{V}) which transform the matrix \textit{M}(\textit{V}) into \textit{M}(\textit{U}).

From the definition of \(h\circ \bar{\sigma }\) it follows that if \(V \subset A^n\) is a neighborly \textit{d}-code, then \(h\circ \bar{\sigma }(V)\) is still a neighborly \textit{d}-code and \(|V|=|h\circ \bar{\sigma }(V)|\). Therefore, in many reasoning we may change an initial code \textit{V} into its isomorphic version \(h\circ \bar{\sigma }(V)\) whose form is more convenient for our purposes than the form of \textit{V}. Below, based on Lemma \hyperlink{FPar3}{1}, we describe such convenient form of \textit{V}.

Let \(V\subset A^n\) be a neighborly code, and let \(j\in [n], \delta \in \{0,1\}\) be such that

\hypertarget{Equ19}{}\begin{equation*}
\begin{aligned} |V^{j,\delta }|=\max _{i\in [n],\varepsilon \in \{0,1\}}|V^{i,\varepsilon }|. \end{aligned}
\end{equation*}

We are intend to work with codes \textit{V} such that \(|V^{j,0}|\ge |V^{i,\varepsilon }|\) for every \(i\in [n]\) and \(\varepsilon \in \{0,1\}\), so if \(|V^{j,0}|<|V^{j,1}|=\max _{i\in [n],\varepsilon \in \{0,1\}}|V^{i,\varepsilon }|\), then we may flip all letters in all words \(v\in V\) at the \textit{j}-th position passing in this way from \textit{V} to its isomorphic form \textit{W} such that \(|W^{j,0}|\ge |W^{i,\varepsilon }|\) for every \(i\in [n]\) and \(\varepsilon \in \{0,1\}\). Due to the possibility of such transition to an isomorphic code, we can immediately assume that the code \textit{V} has the property \(|V^{j,0}|\ge |V^{i,\varepsilon }|\) for \(i\in [n]\) and \(\varepsilon \in \{0,1\}\). Moreover, we assume that \(|V^{j,0}|,|V^{j,1}|\ge 2\). By Lemma \hyperlink{FPar3}{1}, there are disjoint and non-empty sets \(C^j_0,C^j_1\subset [n]{\setminus } \{j\}\) such that \([n]{\setminus } \{j\}=C^j_0\cup C^j_1\cup D\), where \(D=[n]{\setminus } (C^j_0\cup C^j_1\cup \{j\})\) with the following properties:

(\(\alpha \)) For every \(k\in C^j_0\) and \(v\in V^{j,1}\) we have \(v_k=*\). Moreover, for every \(k\in C^j_0\) there are \(u,w\in V^{j,0}\) such that \(u_k+w_k=1\).

(\(\beta \)) For every \(k\in C^j_1\) and \(v\in V^{j,0}\) we have \(v_k=*\). Moreover, for every \(k\in C^j_1\) there are \(u,w\in V^{j,1}\) such that \(u_k+w_k=1\).

(\(\gamma \)) If \(D\ne \emptyset \), then for every \(k\in D\) and for every two words \(v,u\in V^{j,0}\cup V^{j,1}\) we have \(v_k+u_k\ne 1\)

(The property (\(\gamma \)) steams from the definition of the set \textit{D} and the fact that \textit{V} is a neighborly code: If \(v\in V^{j,0}, u \in V^{j,1}\), then, since \(v_j+u_j=1\) and \(k\ne j\), we have \(v_k+u_k\ne 1\)).

For clarity of our notation, again by possibility of passing to an isomorphic form of \textit{V}, we may assume that \textit{V} is such that \(j=1\), that is,

\hypertarget{Equ2}{}\begin{equation*}
\begin{aligned} |V^{1,0}|=\max _{i\in [n],\varepsilon \in \{0,1\}}|V^{i,\varepsilon }|, \end{aligned}\tag{2}
\end{equation*}

and \(C^1_0=\{2, \ldots ,s\}\) for some \(s=s(V)\in [n]\), \(C^1_1=\{s+1, \ldots ,r\}\) for some \(r=r(V)\in [n]\) and finally, if \(D\ne \emptyset \), then \(D=\{r+1, \ldots ,n\}\). Moreover,

(\(\gamma '\)) if \(D\ne \emptyset \), then for every \(k\in D\) and for every \(v\in V^{1,0}\cup V^{1,1}\) we have \(v_k\in \{0,*\}\).

If \textit{V} is as above, then we say that it is in \textit{standard} form (compare Table \hyperlink{Tab1}{1} and the second example in Examples 1).

\hypertarget{FPar5}{}\paragraph{Remark 1.}

Of course, we could work, by Lemma \hyperlink{FPar3}{1}, with codes which are not in standard form but then for example, Table \hyperlink{Tab1}{1} would be far less readable than in the case of codes in standard forms.

We defined standard form for codes with \(C^1_0, C^1_1\ne \emptyset \), which makes the notations easier, as in our proof of Theorem \hyperlink{FPar2}{2} this assumption will be satisfied.

\begin{table}[htbp]\centering
\hypertarget{Tab1}{}\caption{The structure of a neighborly code \textit{V} in standard form in the case \(D\ne \emptyset \), where rows of \textit{M}(\textit{V}) are words in \textit{V}}
\begin{tabular}{cccccccccc}
\hline

1 & 2 & \(\cdots \) & \textit{s} & \(s+1\) & \(\cdots \) & \textit{r} & \(r+1\) & \(\cdots \) & \textit{n} \\\hline

0 & \(v_{2}\) & \(\cdots \) & \(v_{s}\) & \(*\) & \(\cdots \) & \(*\) & \(v_{r+1}\) & \(\cdots \) & \(v_{n}\) \\

\(\vdots \) & \(\vdots \) &  & \(\vdots \) & \(\vdots \) &  & \(\vdots \) & \(\vdots \) &  & \(\vdots \) \\

0 & \(u_{2}\) & \(\cdots \) & \(u_{s}\) & \(*\) & \(\cdots \) & \(*\) & \(u_{r+1}\) & \(\cdots \) & \(u_{n}\) \\

1 & \(*\) & \(\cdots \) & \(*\) & \(w_{s+1}\) & \(\cdots \) & \(w_r\) & \(w_{r+1}\) & \(\cdots \) & \(w_n\) \\

\(\vdots \) & \(\vdots \) &  & \(\vdots \) & \(\vdots \) &  & \(\vdots \) & \(\vdots \) &  & \(\vdots \) \\

1 & \(*\) & \(\cdots \) & \(*\) & \(p_{s+1}\) & \(\cdots \) & \(p_{r}\) & \(p_{r+1}\) & \(\cdots \) & \(p_{n}\) \\

\(*\) & \(t_{2}\) & \(\cdots \) & \(t_{s}\) & \(t_{s+1}\) & \(\cdots \) & \(t_{r}\) & \(t_{r+1}\) & \(\cdots \) & \(t_{n}\) \\

\(\vdots \) & \(\vdots \) &  & \(\vdots \) & \(\vdots \) &  & \(\vdots \) & \(\vdots \) &  & \(\vdots \) \\

\(*\) & \(q_{2}\) & \(\cdots \) & \(q_{s}\) & \(q_{s+1}\) & \(\cdots \) & \(q_{r}\) & \(q_{r+1}\) & \(\cdots \) & \(q_{n}\). \\

\hline\end{tabular}\end{table}

Note that, by the properties (\(\alpha \)),(\(\beta \)) and (\(\gamma '\)), every column in the sub-matrix of \textit{M}(\textit{V}) of the form

\hypertarget{Equ20}{}\begin{equation*}
\begin{aligned}{}[a_{ij}]_{i\in \{1, \ldots ,|V^{1,0}|\},j\in \{2, \ldots ,s\}} \end{aligned}
\end{equation*}

contains at least one 0 and at least one 1, and the sub-matrix

\hypertarget{Equ21}{}\begin{equation*}
\begin{aligned}{}[a_{ij}]_{i\in \{|V^{1,0}|+1, \ldots ,|V^{1,0}|+|V^{1,1}|\},j\in \{2, \ldots ,s\}} \end{aligned}
\end{equation*}

contains only stars. Similarly, the sub-matrix

\hypertarget{Equ22}{}\begin{equation*}
\begin{aligned}{}[a_{ij}]_{i\in \{1, \ldots ,|V^{1,0}|\},j\in \{s+1, \ldots ,r\}} \end{aligned}
\end{equation*}

contains only stars, while every column in the sub-matrix

\hypertarget{Equ23}{}\begin{equation*}
\begin{aligned}{}[a_{ij}]_{i\in \{|V^{1,0}|+1, \ldots ,|V^{1,0}|+|V^{1,1}|\},j\in \{s+1, \ldots ,r\}} \end{aligned}
\end{equation*}

contains at least one 0 and at least one 1. Finally, every column in the sub-matrix

\hypertarget{Equ24}{}\begin{equation*}
\begin{aligned}{}[a_{ij}]_{i\in \{1, \ldots ,|V^{1,0}|+|V^{1,1}|\},j\in \{r+1, \ldots ,n\}} \end{aligned}
\end{equation*}

contains only 0's or stars.

In what follows a flip of a letter \(a\in A=\{0,1,*\}\) will be denoted by \(f(a)=a'\), that is, \(0'=1, 1'=0\) and \(*'=*\).

At the end of this section, we show that in a neighborly \textit{d}-code \(V\subset A^n\) at least one of the sets \(V^{i,a}\), \(i\in [n], a\in A\), is relatively large.

\hypertarget{FPar6}{}\paragraph{Lemma 2.}

If \(V\subset A^n\) is a code, then

\hypertarget{Equ3}{}\begin{equation*}
\begin{aligned} V{\setminus } \{v\}=\bigcup _{j\in {\text {prop}}\, (v)}V^{j,v_j'} \end{aligned}\tag{3}
\end{equation*}

for every \(v\in V\). Consequently, if \textit{V} is a neighborly \textit{d}-code, then for every \(v\in V\) there is \(i\in {\text {prop}}\, (v)\) such that

\hypertarget{Equ4}{}\begin{equation*}
\begin{aligned} |V^{i,v_i'}|\ge (|V|-1)/d. \end{aligned}\tag{4}
\end{equation*}

\hypertarget{FPar7}{}\paragraph{Proof.}

Fix \(v\in V\). Since \(v_j\ne v_j'\) for every \(j\in \textrm{prop}(v)\), we have \(v\not \in V^{j,v_j'}\) for every \(j\in {\text {prop}}\, (v)\). Thus \( V^{j,v_j'}\subset V{\setminus } \{v\}\) for every \(j\in \textrm{prop}(v)\) and therefore \(\bigcup _{j\in {\text {prop}}\, (v)}V^{j,v_j'}\subset V{\setminus } \{v\}\).

To get the opposite inclusion let \(u\in V{\setminus } \{v\}\). Since \textit{V} is a code, there is \(j\in {\text {prop}}\, (v)\) such that \(u_j=v_j'\), that is \(u\in V^{j,v_j'}\). Hence, \(V{\setminus } \{v\}\subset \bigcup _{j\in {\text {prop}}\, (v)}V^{j,v_j'}\), and thus, (\hyperlink{Equ3}{3}) holds true.

To prove (\hyperlink{Equ4}{4}) observe first that it is obvious for \(d=1\), because then we have \(|V|\le 2\).

Let \(d\ge 2\). For every \(i,j\in {\text {prop}}\, (v)\), \(i\ne j\), we have \(V^{i,v_i'}\cap V^{j,v_j'}=\emptyset \). To show this assume on the contrary that \(u\in V^{i,v_i'}\cap V^{j,v_j'}\) for some \(i,j\in {\text {prop}}\, (v)\), \(i\ne j\). Then \(u_i=v_i'\) and \(u_j=v_j'\) which means that the words \textit{u}, \textit{v} are not neighborly which is impossible as \textit{V} is neighborly. Since the sets \((V^{j,v_j'})_{j\in {\text {prop}}\, (v)}\) are pairwise disjoint, we have, by (\hyperlink{Equ3}{3}),

\hypertarget{Equ25}{}\begin{equation*}
\begin{aligned} |V|-1=\sum _{j\in {\text {prop}}\, (v)}|V^{j,v_j'}| \end{aligned}
\end{equation*}

and hence, as \(|{\text {prop}}\, (v)|=d\), we obtain \(|V^{i,v_i'}|\ge (|V|-1)/d\) for some \(i\in {\text {prop}}\, (v)\). \(\square \)

\hypertarget{Sec3}{}\section{An Inflation of a Code}

In this section, we define an inflation of a code which is the main tool in our proof of Theorem \hyperlink{FPar2}{2}.

Let \(V\subset A^n\) be a code, and let \(V=V^{i,0}\cup V^{i,1}\cup V^{i,*}\), where \(i\in [n]\) (recall that \(V^{i,a}=\{v\in V:v_i=a\}\), \(a\in A\)).

For \(\eta =0,1\) let

\hypertarget{Equ26}{}\begin{equation*}
\begin{aligned} V^{i,\eta ,*}=\{v_1\ldots v_{i-1}*v_{i+1}\ldots v_n:v\in V^{i,\eta }\}. \end{aligned}
\end{equation*}

Note that \(V^{i,\eta ,*}\) is a code: For every two words \(v,u \in V^{i,\eta ,*}\) the words \(\bar{u}=u_1\ldots u_{i-1}\eta u_{i+1}\ldots u_n\) and \(\bar{v}=v_1\ldots v_{i-1}\eta v_{i+1}\ldots v_n\) belong to the code \(V^{i,\eta }\) which means that there is \(j\in [n]{\setminus } \{i\}\) such that \(v_j+u_j=1\), \(v_j,u_j\in \{0,1\}\). Thus, \textit{v} and \textit{u} are dichotomous at the \textit{j}-th position. It follows that the sets \(V^{i,0,*}\cup V^{i,*}\) and \(V^{i,1,*}\cup V^{i,*}\) are codes. To show this, let \(v\in V^{i,0,*}\) and \(u\in V^{i,*}\). Then the word \(\bar{v}=v_1 \ldots v_{i-1}0v_{i+1} \ldots v_n\) belongs to \(V^{i,0}\). Since \(V^{i,0}\cup V^{i,*}\) is a code, there is \(j\in [n]{\setminus } \{i\}\) such that the words \(\bar{v},u\) are dichotomous at the \textit{j}-th position, that is, \(\bar{v}_j+u_j=1\). Consequently, \textit{v}, \textit{u} are dichotomous at the \textit{j}-th position as \(\bar{v}_j=v_j\).

\textit{An inflation of V at the i-th position} is the code \(V^{\delta _i}=V^{i,\delta _i,*}\cup V^{i,*}\), where \({\text {vol}}\, (V^{i,\delta _i'})\le {\text {vol}}\, (V^{i,\delta _i})\), \(\delta _i,\delta _i'\in \{0,1\}\) and \(\delta _i'+\delta _i=1\).

Let \(J=\{i_1< \cdots <i_m\}\subset [n]\), and let the sequence \(J_1=(j_1, \ldots ,j_m)\) be a permutation of elements of the set \textit{J}.

The code

\hypertarget{Equ27}{}\begin{equation*}
\begin{aligned} V^\delta =( \ldots (V^{\delta _{j_1}})^{\delta _{j_2}} \ldots )^{\delta _{j_m}}, \end{aligned}
\end{equation*}

where \(\delta =(\delta _{j_k})_{k\in [m]}\in \{0,1\}^m\), is called an \textit{inflation of V} on the sequence \(J_1\). The sequence \(\delta \) is called an \textit{inflation sequence}. By the definition of inflation, we have

\hypertarget{Equ28}{}\begin{equation*}
\begin{aligned} {\text {vol}}\, (V)\le {\text {vol}}\, (V^\delta ) \end{aligned}
\end{equation*}

for every sequence \(J_1=(j_1, \ldots ,j_m)\) and every inflation sequence \(\delta =(\delta _{j_k})_{k\in [m]}\). Usually we shall indicate only the set \textit{J} without specifying a permutation of \textit{J}. In such a case we just say that an inflation of \textit{V} is on the set \textit{J}. However, as we show in the second part of Examples 1 an inflation of a code depends on a permutation of elements of \textit{J}.

Let \(\bar{J}=\{i_1< \cdots <i_r\}\subset J\), \(\bar{J}_1=(j_1, \ldots ,j_r)\), where \(r<m\), and let \(\bar{\delta }=(\delta _{j_k})_{k\in [r]}\) be an inflation sequence on the sequence \(\bar{J}_1\). Fix any \(i\in J{\setminus } \bar{J}\) and let us consider an inflation on the sequence \(\tilde{J}_1=(j_1, \ldots ,j_r,i)\). Observe that the code \(V^{\bar{\delta }}\) for such \textit{i} can be in one of the three states:

\hypertarget{Equ29}{}\begin{equation*}
\begin{aligned}{} & {} {\text {vol}}\, ((V^{\bar{\delta }})^{i,0})>{\text {vol}}\, ((V^{\bar{\delta }})^{i,1})\,\, \textrm{or}\,\, \\{} & {} \quad {\text {vol}}\, ((V^{\bar{\delta }})^{i,0})<{\text {vol}}\, ((V^{\bar{\delta }})^{i,1}) \,\, \textrm{or}\,\, {\text {vol}}\, ((V^{\bar{\delta }})^{i,0})={\text {vol}}\, ((V^{\bar{\delta }})^{i,1}). \end{aligned}
\end{equation*}

In the first two cases we say that \(V^{\bar{\delta }}\) is in \textit{0-advantage} (resp. \textit{1-advantage}) at the \textit{i}-th position. In the third case we say that the code \(V^{\bar{\delta }}\) is \textit{balanced} at the \textit{i}-th position.

Let \(\tilde{\delta }=(\delta _{j_1}, \ldots ,\delta _{j_k},\delta _i)\). If \(V^{\bar{\delta }}\) is in 0-advantage at the \textit{i}-th position, then \(\delta _i=0\), and consequently, all words from the set \((V^{\bar{\delta }})^{i,1}\) have to be removed, and all words from the set \((V^{\bar{\delta }})^{i,0}\) have to be modified. This means that in the code \(V^{\bar{\delta }}=(V^{\bar{\delta }})^{i,0}\cup (V^{\bar{\delta }})^{i,1}\cup (V^{\bar{\delta }})^{i,*}\) the set of words \((V^{\bar{\delta }})^{i,1}\) is removed, and every word in \((V^{\bar{\delta }})^{i,0}\) is modified by changing every 0 to \(*\) at the \textit{i}-th position. In the result we obtain the inflation on the sequence \(\tilde{J}_1\) which is of the form \(V^{\tilde{\delta }}=(V^{\bar{\delta }})^{i,0,*}\cup (V^{\bar{\delta }})^{i,*}\). If \(V^{\bar{\delta }}\) is balanced at the \textit{i}-th position, then we have a choice: We may take \(\delta _i=0\) or \(\delta _i=1\). In this case we get \(V^{\tilde{\delta }}=(V^{\bar{\delta }})^{i,0,*}\cup (V^{\bar{\delta }})^{i,*}\) (if \(\delta _i=0\)) or \(V^{\tilde{\delta }}=(V^{\bar{\delta }})^{i,1,*}\cup (V^{\bar{\delta }})^{i,*}\) (if \(\delta _i=1\)).

Thus, any inflation \(V^\delta \) of a code \textit{V} on a set \(J\subset [n]\) is a code that arises from \textit{V} in such a way that some words of \textit{V} are removed, some are modified and some words from the code \textit{V} are unmodified. Therefore, for every \(u\in V^\delta \) there is \(v\in V\) such that \(u_i=v_i\) for every \(i\in {\text {prop}}\, (u)\subset {\text {prop}}\, (v)\). If \({\text {prop}}\, (u)\subsetneq {\text {prop}}\, (v)\), then we say that \textit{u} is a \textit{modification} of the word \textit{v}. If \({\text {prop}}\, (u)={\text {prop}}\, (v)\) then we say that \textit{v} is \textit{unmodified} during an inflation process on \textit{J}. In this case, by the definition of inflation, \(v_i=*\) for every \(i\in J\) (compare the second example below).

\hypertarget{FPar8}{}\paragraph{Example 1.}

Let \(V=\{00*,*11\}\) (Fig. \hyperlink{Fig1}{1}). The code \textit{V} is balanced at the position 2, it is in 1-advantage state at the position 3 and in 0-advantage state at the position 1. Let \(J=\{2\}\). Then \(\delta =(\delta _i)_{i\in J}\) is an inflation sequence, where \(\delta _2=0\) and \(V^\delta =\{0**\}\). Of course, if \(\delta _2=1\), then \(\delta \) is also an inflation sequence and \(V^\delta =\{**1\}\). In both cases we have \({\text {vol}}\, (V^\delta )={\text {vol}}\, (V)\). Let \(J=\{3\}\). Then \(\delta =(\delta _i)_{i\in J}\) is an inflation sequence for \(\delta _3=1\) and \(V^\delta =\{00*,*1*\}\). In this case \({\text {vol}}\, (V^{\delta })>{\text {vol}}\, (V)\). If \(J=\{3\}\) and \(\delta _3=0\), then \(\delta \) is not an inflation sequence as \(0={\text {vol}}\, (V^{3,0})<{\text {vol}}\, (V^{3,1})=2\).


\begin{figure}[htbp]\centering
\hypertarget{Fig1}{}\resizebox{0.94\linewidth}{!}{\InflationFigure}
\caption{Let \(V=\{00*, *11\}\). Then \(\breve{V}=\{[0,1]^2\times [0,2],[0,2]\times [1,2]^2]\) (the picture on the left). We have \(V^\delta =\{0**\}\) for \(J=\{2\}\) and \(\delta _2=0\) (a realization of \(V^\delta \) is given in the middle) and \(V^\delta =\{00*, *1*\}\) for \(J=\{3\}\) and \(\delta _3=1\) (a realization of \(V^\delta \) is given on the right)}\end{figure}

Our second example concerns the following neighborly code \(W\subset A^6\) (note that, \textit{W} is in standard form):

\begin{table}[htbp]\centering
\hypertarget{Tab2}{}\caption{A neighborly code in standard form, where \(C^1_0=\{2,3,4\}, C^1_1=\{5\}\) and \(D=\{6\}\)}
\begin{tabular}{cccccc}
\hline

1 & 2 & 3 & 4 & 5 & 6 \\\hline

0 & \(*\) & 1 & \(*\) & \(*\) & 0 \\

0 & 1 & 0 & \(*\) & \(*\) & 0 \\

0 & 0 & 0 & 1 & \(*\) & \(*\) \\

0 & 0 & 0 & 0 & \(*\) & 0 \\

1 & \(*\) & \(*\) & \(*\) & 0 & 0 \\

1 & \(*\) & \(*\) & \(*\) & 1 & \(*\) \\

\(*\) & \(*\) & \(*\) & 0 & 0 & 1 \\

\hline\end{tabular}\end{table}

An inflation of a code usually depends on a sequence on which it is made, that is, for a given sequence \(J_1\) if \(J_2\) is a permutation of \(J_1\), then it can happen that inflation on \(J_1\) is not equal to the inflation on \(J_2\). For example, for the code \textit{W} given in Table \hyperlink{Tab2}{2} we let \(J=\{1,2,3\}\), \(J_1=(1,2,3)\) and \(J_2=(3,2,1)\). The inflation of \textit{W} on the sequence \(J_1\) is of the form \(W^{\delta ^1}=\{****00,****1*,***001\}\), where \(\delta ^1=(1,0,0)\), while the inflation of \textit{W} on the sequence \(J_2\) is of the form \(W^{\delta ^2}=\{***1**,***0*0,***001\}\), where \(\delta ^2=(0,0,0)\). Thus, \(W^{\delta ^1}\ne W^{\delta ^2}\). (Note that, we may take \(\delta ^3=(0,0,1)\) on \(J_2\), as in the last step we have a balance at the 1-th position, and then \(W^{\delta ^1}=W^{\delta ^3}\).)

During the inflation process of \textit{W} on \(J_2\) the word \(v=010**0\) is (in the first step) modified to the word \(u=01***0\), but in the second step of the inflation process on \(J_2\), the word \textit{u} is removed (in this sense \textit{v} is removed during an inflation process). On the other hand, the word \(w=***001\) is unmodified during the inflation process on \(J_1\) and \(J_2\). The first two words in \(W^{\delta ^2}\) are modifications of the third and the fourth word in \textit{W}, respectively, and the first two words in \textit{W} as well as the fifth and sixth words in \textit{W} are removed during the inflation process on \(J_2\) for the inflation sequence \(\delta ^2\).

To conclude this section, we discuss briefly the behavior of a family of neighborly simplices when corresponding \textit{d}-code is subjected to an inflation process. Figure \hyperlink{Fig2}{2} presents three pairwise neighborly triangles \(T_1,T_2\) and \(T_3\), where 0's and 1's denote the sides of lines \(H_1, \ldots ,H_6\). (Below Fig. \hyperlink{Fig2}{2} these sides are denoted by \(H_i^{\varepsilon }\), where \(\varepsilon \in \{0,1\}\).) Let \(v^1=v(T_1), v^2=v(T_2), v^3=v(T_3)\) be words as described at the beginning of Sect. \hyperlink{Sec1}{1}. We have

\hypertarget{Equ30}{}\begin{equation*}
\begin{aligned}{} & {} v^1=0**10*, \\{} & {} v^2=*110**, \\{} & {} v^3=***111. \end{aligned}
\end{equation*}

Let \(J=\{4\}\) and \(\delta =(1)\). The inflation \(V^\delta \) of the code \(V=\{v^1,v^2,v^3\}\) on the set \textit{J} with the inflation sequence \(\delta \) is of the form \(V^\delta =\{0***0*,****11\}\). Thus, the word \(v^2\) is removed and two words are modified. For triangles this inflation process means that the triangle \(T_2\) is removed and instead of triangles \(T_1\) and \(T_3\) we obtain two cons \(C_1\) and \(C_3\) which correspond to words from \(V^\delta \) (as presented in Fig. \hyperlink{Fig2}{2}). The usefulness of an inflation of a code in examinations of the structure of codes lies on two things: The first is the inequality \({\text {vol}}\, (V)\le {\text {vol}}\, (V^\delta )\) and the second is the ability to control of volumes of boxes at individual stages of the inflation process. In the case of triangles \(T_1,T_2,T_3\) we lost information on areas, as the areas of cones \(C_1,C_2\) are infinite. Moreover, in the case of codes during the whole inflation process we work with neighborly boxes, while in the case of a family of neighborly simplices some of the simplices of the family are turned into unbounded convex polyhedrons. On the other hand, such modified family is still neighborly and perhaps for this reason an inflation of a code can be useful in examinations of neighborly families of simplices.

\begin{figure}[htbp]\centering
\hypertarget{Fig2}{}\resizebox{0.65\linewidth}{!}{\SimplicesFigure}
\caption{Simplices \(T_1=H_1^0\cap H_4^1\cap H_5^0, T_2=H_2^1\cap H_3^1\cap H_4^0, T_3=H_4^1\cap H_5^1\cap H_6^1\) and cones \(C_1=H_1^0\cap H_5^0\) and \(C_3=H_5^1\cap H_6^1\)}\end{figure}

\hypertarget{Sec4}{}\section{A Proof of Theorem 2}

Our proof of Theorem \hyperlink{FPar2}{2} consists in controlling some inflation process of a neighborly \textit{d}-code \textit{V} without twin pairs in such a way that some portion of \textit{V} remains unmodified (it will be \(V^{1,1}\)) during the inflation process, and on the other hand the form of some portion of the considered inflation of \textit{V} is easy to predict.

\hypertarget{FPar9}{}\paragraph{Proof of Theorem 2.}

Suppose that the theorem is not true. Then there are an integer \(M>0\), a sequence of positive integers \((d_m)_{m\ge 1}\) with \(d_m\rightarrow \infty \) as \(m\rightarrow \infty \) and a sequence \(V_{d_m}\subset A^{n(d_m)}\) of neighborly \(d_m\)-codes without twin pairs such that \(|V_{d_m}|=2^{d_m}-M\) for \(m\ge 1\). We may assume that \textit{M} is the smallest such number, that is, there is \(d_0\ge 1\) such that for every \(d\ge d_0\) and every neighborly \textit{d}-code \(W\subset A^{n(d)}\) if \(|W|=2^d-(M-k)\), where \(k\in [M]\), then \textit{W} contains a twin pair. By (\hyperlink{Equ1}{1}) we have \(M\ge 2\).

Let \(d\in \{d_m:m\in \mathbb {N}\}\) be such that \(d>d_0\), and let \(V=V_{d}\subset A^{n(d)}\) be a neighborly \textit{d}-code without twin pairs with \(|V|=2^{d}-M\).

As we show below, \(|V^{1,0}|,|V^{1,1}|\ge 2\), and thus we may assume that \textit{V} is in standard form (compare Sect. \hyperlink{Sec2}{2} and Table \hyperlink{Tab1}{1} and \hyperlink{Tab2}{2}).

By (\hyperlink{Equ2}{2}) and (\hyperlink{Equ4}{4}), we have \(|V^{1,0}|\ge (2^d-(M+1))/d\). Therefore, (increasing \(d\in \{d_m:m\in \mathbb {N}\}\), where \(d>d_0\), if needed) we may assume that

\hypertarget{Equ5}{}\begin{equation*}
\begin{aligned} |V^{1,0}|> 9M. \end{aligned}\tag{5}
\end{equation*}

Since \(|V|=2^d-M\), we have

\hypertarget{Equ6}{}\begin{equation*}
\begin{aligned} |V^{1,1}|\ge |V^{1,0}|-M. \end{aligned}\tag{6}
\end{equation*}

To show this, let \(W^{1,1}=\{1v_2 \ldots v_n:v\in V^{1,0}\}\). The set \(W^{1,1}\) is a code because \(V^{1,0}\) is a code, and since the set \(V^{1,0}\cup V^{1,*}\) is a code, the set \(W^{1,1}\cup V^{1,*}\) must be a code. Therefore, the set \(W=V^{1,0}\cup W^{1,1}\cup V^{1,*}\) is a code. Moreover, since \textit{V} is a \textit{d}-code, \textit{W} is a \textit{d}-code, and thus \(|W|\le 2^d\) (see Introduction). If on the contrary \(|V^{1,1}|< |V^{1,0}|-M\), then since \(|W^{1,1}|=|V^{1,0}|\) and \(|V|=2^{d}-M\), we have

\hypertarget{Equ31}{}\begin{equation*}
\begin{aligned} |W|=|V^{1,0}|+|V^{1,0}|+|V^{1,*}|>|V^{1,0}|+|V^{1,1}|+|V^{1,*}|+M =|V|+M=2^d, \end{aligned}
\end{equation*}

which is a contradiction.

In what follows we shall consider an inflation \(V^\delta =U\) on the set \(J=C^1_0=\{2,\ldots ,s\}\). Note that, by the property (\(\alpha \)) in Sect. \hyperlink{Sec2}{2}, for every \(v\in V^{1,1}\) and every \(k\in J\) we have \(v_k=*\). Thus, each word in \(V^{1,1}\) is unmodified and of course is not removed during an inflation process on the set \textit{J}. Therefore,

\hypertarget{Equ32}{}\begin{equation*}
\begin{aligned} U^{1,1}=V^{1,1} \end{aligned}
\end{equation*}

(compare Table \hyperlink{Tab1}{1} and Examples 1). Moreover, since \textit{V} is in standard form, and the inflation of \textit{V} is on the set \textit{J}, it follows that \(|U^{1,0}|\le 1\). Indeed, suppose on the contrary that there are two words \(v,u\in U^{1,0}\). The words \textit{v}, \textit{u} arose during the inflation process from some two words \(p,q\in V^{1,0}\), but modifications of \textit{p}, \textit{q} to \textit{v}, \textit{u} were made only on the set \textit{J}. This means, taking into account the property (\(\beta \)) given in Sect. \hyperlink{Sec2}{2}, that \(v_i=u_i=*\) for every \(i\in \{2,\ldots ,r\}\) (as, by (\(\beta \)), \(p_k=q_k=*\) for \(k\in \{s+1,\ldots ,r\}\)). The set \textit{U} is a code, and therefore \textit{v}, \textit{u} are dichotomous, that is \(v_i+u_i=1\), where \(i\in \{r+1,\ldots ,n\}\). Since \(v_i=p_i\) and \(q_i=u_i\), we obtain \(p_i+q_i=1\) for \(i\in \{r+1,\ldots ,n\}\). This contradicts the property (\(\gamma '\)) given in Sect. \hyperlink{Sec2}{2}. Therefore, \(U^{1,0}=\emptyset \) or \(U^{1,0}\) contains precisely one word.

Now we consider three cases depending on the form of \(U^{1,0}\).

\textbf{Case 1}. Let us suppose that there is an inflation \(V^\delta =U\) on the set \textit{J} such that

\hypertarget{Equ33}{}\begin{equation*}
\begin{aligned} U^{1,0}=\emptyset , \end{aligned}
\end{equation*}

that is, \(U=V^{1,1}\cup U^{1,*}\). Since \({\text {vol}}\, (U)\ge {\text {vol}}\, (V)\), we have \({\text {vol}}\, (U)\ge (2^d-M)2^{n-d}\). Let \(W=W^{1,0}\cup V^{1,1}\cup U^{1,*}\), where \(W^{1,0}=\{0v_2\ldots v_n:v\in V^{1,1}\}\). Clearly, \(W\subset A^n\) is a code (we show this in the similar manner as in the case of the code \textit{W} given right after (\hyperlink{Equ6}{6})) and hence \({\text {vol}}\, (W)\le 2^n\). Moreover, \({\text {vol}}\, (W^{1,0})={\text {vol}}\, (V^{1,1})\), by the definition of \(W^{1,0}\). On the other hand, since, by (\hyperlink{Equ5}{5}) and (\hyperlink{Equ6}{6}), \(|V^{1,1}|>M\) and \(|v|=2^{n-d}\) for \(v\in V^{1,1}\) (as \(V^{1,1}\) is a \textit{d}-code), we have \({\text {vol}}\, (V^{1,1})>M2^{n-d}\), and consequently

\hypertarget{Equ34}{}\begin{equation*}
\begin{aligned} {\text {vol}}\, (W)={\text {vol}}\, (U)+{\text {vol}}\, (V^{1,1})>(2^d-M)2^{n-d}+M2^{n-d}=2^n, \end{aligned}
\end{equation*}

a contradiction.

\textbf{Case 2}. We now assume that for each inflation \(V^\delta =U\) on the set \textit{J} we have

\hypertarget{Equ35}{}\begin{equation*}
\begin{aligned} U^{1,0}=\{0*\ldots *\}. \end{aligned}
\end{equation*}

Then \(U^{1,*}=\emptyset \) as \textit{U} is a code (if \(v\in U^{1,*}\), then the words \textit{v} and \(0*\ldots *\) are not dichotomous because \(v_1=*\)). Thus, \(U=\{0*\ldots *\}\cup V^{1,1}\). We shall show that this form of \textit{U} is not possible.

Let us assume first that there is at least one inflation \(V^\delta =U\) such that \({\text {vol}}\, (U)>{\text {vol}}\, (V)\) (in other words, during the inflation process defined by \(\delta \) there are unbalanced states). Since \({\text {vol}}\, (U)>{\text {vol}}\, (V)=2^n-M2^{n-d}\) and \({\text {vol}}\, (U^{1,0})={\text {vol}}\, (\{0*\ldots *\})=2^{n-1}\), we have

\hypertarget{Equ36}{}\begin{equation*}
\begin{aligned} {\text {vol}}\, (V^{1,1})>2^{n-1}-M2^{n-d}. \end{aligned}
\end{equation*}

The set \(W=\{v_2\ldots v_n:v\in V^{1,1}\}\) is a neighborly code with \(|{\text {prop}}\, (v)|=d-1\) for every \(v\in W\) and \({\text {vol}}\, (W)={\text {vol}}\, (V^{1,1})\) as \(v_1=1\) for every \(v\in V^{1,1}\). Thus, \({\text {vol}}\, (W)>2^{n-1}-M2^{n-1-(d-1)}\), and since \({\text {vol}}\, (W)=|W|2^{n-1-(d-1)}\), we obtain \(|W|=2^{d-1}-(M-k)\) for some \(k\in [M]\). Since \(d-1\ge d_0\) (recall that \(d>d_0\)), by the minimality of \textit{M}, the code \textit{W} contains a twin pair, say \textit{v}, \textit{u}, and then \(V^{1,1}\) contains the twin pair 1\textit{v}, 1\textit{u}, where \(1v=1v_2\ldots v_n\) and \(1u=1u_2\ldots u_n\). A contradiction.

Therefore, since now we assume that for each inflation \(V^\delta =U\) on \textit{J} we have \({\text {vol}}\, (U)={\text {vol}}\, (V)\) (that is, any \(0-1\) sequence \(\delta \) is an inflation sequence).

Let us suppose that \(V^{1,*}\ne \emptyset \). Take any \(q\in V^{1,*}\), and let \( B\subset J\) be such that \(q_j\in \{0,1\}\) for \(j\in B\). Note that, \(B\ne \emptyset \) as if \(B=\emptyset \), that is \(q=*\ldots *q_{s+1}\ldots q_n\), then \(q\in U^{1,*}\) (such word \textit{q} cannot be removed during any inflation process on the set \(J=\{2,\ldots ,s\}\) as \(q_j=*\) for all \(j\in J\)). Now we consider an inflation on the sequence \(J_1=(2,\ldots ,s)\) with \(\delta =(\delta _j)_{j\in J}\) such that \(\delta _j=q_j\) for \(j\in B\). By our assumption, \(\delta \) is an inflation sequence, and thus \(*\ldots *q_{s+1}\ldots q_n\in U^{1,*}\), which contradictions the equality \(U^{1,*}=\emptyset \).

Let now \(V^{1,*}=\emptyset \), that is, \(V=V^{1,0}\cup V^{1,1}\). Since \({\text {vol}}\, (U)=2^{n-1}+{\text {vol}}\, (V^{1,1})\), \({\text {vol}}\, (V)={\text {vol}}\, (V^{1,0})+{\text {vol}}\, (V^{1,1})\) and \({\text {vol}}\, (U)={\text {vol}}\, (V)\), it follows that \({\text {vol}}\, (V^{1,0})=2^{n-1}\). Thus, the neighborly \((d-1)\)-code \(W=\{v_2\ldots v_n:v\in V^{1,0}\}\) contains \(2^{d-1}\) words, and therefore (we have \(M\ge 2\) by (\hyperlink{Equ1}{1})) it contains a twin pair. Hence \(V^{1,0}\) contains a twin pair which is a contradiction.

\textbf{Case 3}. In the last case we assume that there is an inflation \(V^\delta =U\) on \textit{J} such that

\hypertarget{Equ37}{}\begin{equation*}
\begin{aligned} U^{1,0}=\{0*\ldots *u_{r+1}\ldots u_{i-1}0u_{i+1}\ldots u_n\}. \end{aligned}
\end{equation*}

Let us note that the form of the set \(U^{1,0}\) means that \(D\ne \emptyset \), as if \(D=\emptyset \), then \(U^{1,0}\) would be as in Case 1 or 2. Thus, by the property (\(\gamma '\)), \(u_k\in \{0,*\}\) for \(k\in \{r+1,\ldots ,n\}\). As in the previous two cases, we show that the form of \(U^{1,0}\) is not possible. As always \(U^{1,1}=V^{1,1}\).

We shall show first that

\hypertarget{Equ7}{}\begin{equation*}
\begin{aligned} |V^{i,*}\cap V^{1,1}|\le 2M. \end{aligned}\tag{7}
\end{equation*}

We assume that \(V^{i,*}\cap V^{1,1}\ne \emptyset \), otherwise there is nothing to prove. Let \(W=\{0v_2\ldots v_{i-1}1v_{i+1}\ldots v_n:v\in V^{i,*}\cap V^{1,1}\}\). It is easy to see that \textit{W} is a code: If we take two words \textit{v}, \textit{u} belonging to the code \(V^{i,*}\cap V^{1,1}\) then, since \(v_1=u_1=1\) and \(v_i=u_i=*\), we have \(v_k+u_k=1\) for some \(k\in [n]{\setminus }\{ 1,i\}\). Moreover, \({\text {vol}}\, (W)=|V^{i,*}\cap V^{1,1}|2^{n-(d+1)}\) as \(|w|=|v|/2\) for \(w\in W\) and \(v\in V^{i,*}\cap V^{1,1}\) (equivalently, \textit{W} is a \((d+1)\)-code, by the definition of \textit{W}).

Now we show that \(W\cup U\) is a code. Let \textit{w}, \textit{u} be two words such that \(w\in W\) and \(u\in U\). If \(u\in U^{1,0}\), then \(u_i=0\), and since \(w_i=1\), the words \textit{w}, \textit{u} are dichotomous at \textit{i}. If \(u\in U^{1,1}\), then \(u_1=1\). But \(w_1=0\), and then \textit{w}, \textit{u} are dichotomous at the position 1. Finally, let \(u\in U^{1,*}\). The set \textit{U} is a code and \(U^{1,1}=V^{1,1}\) and therefore, if \(v\in V^{1,1}\) is such that \(v_i=*\), then \(u_k+v_k=1\) for some \(k\ne 1,i\). It follows that \textit{w} and \textit{u} are dichotomous at the position \(k\ne 1,i\). Thus, \(W\cup U\) is a code, and hence \({\text {vol}}\, (W\cup U)={\text {vol}}\, (W)+{\text {vol}}\, (U)\). Therefore, \({\text {vol}}\, (W)\le M2^{n-d}\) as \({\text {vol}}\, (U)\ge {\text {vol}}\, (V)=2^n-M2^{n-d}\) and \({\text {vol}}\, (W\cup U)\le 2^n\). Hence,

\hypertarget{Equ38}{}\begin{equation*}
\begin{aligned} |V^{i,*}\cap V^{1,1}|2^{n-(d+1)}\le M2^{n-d} \end{aligned}
\end{equation*}

which gives (\hyperlink{Equ7}{7}).

By the property (\(\gamma '\)) we have \(|V^{1,1}|=|V^{i,*}\cap V^{1,1}|+|V^{i,0}\cap V^{1,1}|\), and thus, from (\hyperlink{Equ7}{7}) and (\hyperlink{Equ6}{6}), it follows that

\hypertarget{Equ8}{}\begin{equation*}
\begin{aligned} |V^{i,0}\cap V^{1,1}|\ge |V^{1,0}|-3M. \end{aligned}\tag{8}
\end{equation*}

By the maximality of \(|V^{1,0}|\) (\hyperlink{Equ2}{2}) we have \(|V^{1,0}|\ge |V^{i,0}|\ge |V^{i,0}\cap V^{1,0}|+|V^{i,0}\cap V^{1,1}|\) and from (\hyperlink{Equ8}{8}) we get

\hypertarget{Equ9}{}\begin{equation*}
\begin{aligned} |V^{i,0}\cap V^{1,0}|\le 3M. \end{aligned}\tag{9}
\end{equation*}

Finally, by the equality \(|V^{1,0}|=|V^{i,*}\cap V^{1,0}|+|V^{i,0}\cap V^{1,0}|\) (which stems from the property (\(\gamma '\))) and from (\hyperlink{Equ9}{9}) we obtain

\hypertarget{Equ10}{}\begin{equation*}
\begin{aligned} |V^{i,*}\cap V^{1,0}|\ge |V^{1,0}|-3M. \end{aligned}\tag{10}
\end{equation*}

Let \(P=\{1v_2\ldots v_{i-1}1v_{i+1}\ldots v_n:v\in V^{i,*}\cap V^{1,0}\}\), and let

\hypertarget{Equ39}{}\begin{equation*}
\begin{aligned} \bigcup \breve{P}=\bigcup _{v\in P}\breve{v}\quad \textrm{and} \quad \bigcup \breve{Q}=\bigcup _{v\in V^{i,*}\cap V^{1,1}}\breve{v}. \end{aligned}
\end{equation*}

(Let us recall that if \(v=v_1\ldots v_n\in A^n\), then \(\breve{v}=\breve{v}_1\times \cdots \times \breve{v}_n\subset [0,2]^n\), where \(\breve{v}_i=[0,1)\) if \(v_i=0\), \(\breve{v}_i=[1,2]\) if \(v_i=1\) and \(\breve{v}_i=[0,2]\) if \(v_i=*\).) Since \textit{V} is a \textit{d}-code, by the definition of \textit{P}, the set \textit{P} is a \((d+1)\)-code, and then

\hypertarget{Equ11}{}\begin{equation*}
\begin{aligned} {\text {vol}}\, (P)=|V^{i,*}\cap V^{1,0}|2^{n-(d+1)}. \end{aligned}\tag{11}
\end{equation*}

Moreover, for every \(w\in P\) and \(v\in V{\setminus }(V^{i,*}\cap V^{1,1})\) the words \textit{w} and \textit{v} are dichotomous. To show this fix \(w=1w_2\ldots w_{i-1}1w_{i+1}\ldots w_n\in P\). If \(v\in V^{1,0}\), then \textit{v}, \textit{w} are dichotomous at the first position; if \(v\in V^{1,*}\), then, by the definition of \textit{P}, \(\bar{w}=0w_2\ldots w_{i-1}*w_{i+1}\ldots w_n\in V^{1,0}\), and therefore there is \(k\in [n]{\setminus } \{1,i\}\) such that \(v_k+w_k=1\), that is, \textit{v}, \textit{w} are dichotomous at the \textit{k}-th position. Finally, if \(v\in V^{1,1}\) and \(v_i=0\) (recall that, by (\(\gamma '\)), \(v_i\in \{0,*\}\)), then \textit{v} and \textit{w} are dichotomous at the \textit{i}-th position. It follows from the above that \(\breve{w}\cap \breve{v}=\emptyset \) for every \(w\in P\) and \(v\in V{\setminus }(V^{i,*}\cap V^{1,1})\). Consequently

\hypertarget{Equ40}{}\begin{equation*}
\begin{aligned} \bigcup \breve{P}{\setminus } \bigcup \breve{Q}\subset [0,2]^n{\setminus } \bigcup _{v\in V}\breve{v}. \end{aligned}
\end{equation*}

Thus, as \(m_n(\bigcup _{v\in V}\breve{v})=2^n-M2^{n-d}\), we obtain

\hypertarget{Equ12}{}\begin{equation*}
\begin{aligned} m_n(\bigcup \breve{P}{\setminus } \bigcup \breve{Q})\le M2^{n-d}. \end{aligned}\tag{12}
\end{equation*}

Since \(m_n(\bigcup \breve{P})={\text {vol}}\, (P)\) and \(m_n(\bigcup \breve{Q})={\text {vol}}\, (V^{i,*}\cap V^{1,1})\), by (\hyperlink{Equ5}{5}), (\hyperlink{Equ7}{7}), (\hyperlink{Equ10}{10}), (\hyperlink{Equ11}{11}) and (\hyperlink{Equ12}{12}) we have

\hypertarget{Equ41}{}\begin{equation*}
\begin{aligned} M2^{n-d}\ge & {} m_n(\bigcup \breve{P}{\setminus } \bigcup \breve{Q})\ge {\text {vol}}\, (P)\ -{\text {vol}}\, (V^{i,*}\cap V^{1,1})\\\ge & {} |V^{i,*}\cap V^{1,0}|2^{n-(d+1)}-2M2^{n-d} \\\ge & {} (|V^{1,0}|-3M)2^{n-(d+1)}-2M2^{n-d}\\= & {} \frac{1}{2}|V^{1,0}|2^{n-d}-\frac{3}{2}M2^{n-d}-2M2^{n-d}>M2^{n-d}, \end{aligned}
\end{equation*}

a contradiction. The proof is completed. \(\square \)

Now our proof of Theorem \hyperlink{FPar1}{1} is immediate:

\hypertarget{FPar10}{}\paragraph{Proof of Theorem 1.}

Since \(S_d\le M_{d+1}\) for every \(d\ge 2\), we have \(2^{d+1}-S_d\ge 2^{d+1}-M_{d+1}\) for \(d\ge 2\). By Theorem \hyperlink{FPar2}{2}, we have \(\lim _{d\rightarrow \infty }(2^{d+1}-M_{d+1})=\infty \), and therefore \(\lim _{d\rightarrow \infty }(2^{d+1}-S_{d})=\infty \) \(\square \)

\hypertarget{Sec5}{}\section{A Result on the Structure of Neighborly \textit{d}-Codes}

From the proof of Theorem \hyperlink{FPar2}{2} we obtain an interesting result on the structure of neighborly \textit{d}-codes:

\hypertarget{FPar11}{}\paragraph{Corollary 1.}

Let \(V\) be a neighborly \(d\)-code without twin pairs in standard form, and let \(M\ge 2\) be an integer such that \(|V|=2^d-M\) and \(|V^{1,0}|> 9\,M\). Then for every inflation \(V^\delta \) on the set \(J=C^1_0=\{2,\ldots ,s\}\) we have

\hypertarget{Equ42}{}\begin{equation*}
\begin{aligned} V^{\delta }=\{0*\ldots *\}\cup V^{1,1}. \end{aligned}
\end{equation*}

Moreover, if \(V^{1,*}\ne \emptyset \), then \({\text {vol}}\, (V^\delta )>{\text {vol}}\, (V)\) for every inflation \(V^\delta \), and consequently \(|V^{1,*}|<M\) and \(|V^{1,0}|, |V^{1,1}|>2^{d-1}-M\).

\hypertarget{FPar12}{}\paragraph{Proof.}

By the proof of Theorem \hyperlink{FPar2}{2}, \(V^\delta =\{0*\ldots *\}\cup V^{1,1}\) for every inflation \(V^\delta \) on \textit{J}. Moreover, if \(V^{1,*}\ne \emptyset \), then, again from the proof of Theorem \hyperlink{FPar2}{2}, \({\text {vol}}\, (V^\delta )>{\text {vol}}\, (V)\) for at least one inflation \(V^\delta \), and by the form of \(V^\delta \), we have \({\text {vol}}\, (V^\delta )>{\text {vol}}\, (V)\) for each inflation \(V^\delta \) on \textit{J}.

Since, \({\text {vol}}\, (\{0*\ldots *\}\cup V^{1,1})>2^n-M2^{n-d}\), that is, \(2^{n-1}+|V^{1,1}|2^{n-d}>2^n-M2^{n-d}\), we have \(|V^{1,1}|>2^{d-1}-M\). But \(|V^{1,0}|\ge |V^{1,1,}|\), and therefore \(|V^{1,0}|>2^{d-1}-M\). Thus, \(|V^{1,*}|<M\) as \(|V^{1,0}|+|V^{1,1}|+|V^{1,*}|=2^d-M\). \(\square \)

\hypertarget{FPar13}{}\paragraph{Remarks 2.}

In [\hyperlink{ref-CR9}{9}] we conjecture that \(M_d=3\cdot 2^{d-2}\). (It is verified for \(d\le 3\).) We believe additionally that, up to isomorphism, there is only one neighborly \textit{d}-code without twin pairs \textit{V} such that \(|V|=3\cdot 2^{d-2}\) and \(V=V^{1,0}\cup V^{1,1}\) (its form is given in [\hyperlink{ref-CR9}{9}, Section 3]). Thus, if it is true that \(M_d=3\cdot 2^{d-2}\), then \(S_d\le M_{d+1}=2^d+2^{d-1}\), This is still far away to \(S_d=2^d\). On the hand, a deeper examination of the structure of neighborly \textit{d}-codes without twin pairs may help to get a better estimation of \(S_d\) than \(S_d\le 2^d+2^{d-1}\) (compare [\hyperlink{ref-CR8}{8}]), if, of course, the equality \(M_d=3\cdot 2^{d-2}\) holds. For example, let us suppose that \(M_d=3\cdot 2^{d-2}\) and the only one, up to isomorphism, neighborly \textit{d}-code without twin pairs \textit{V} with \(|V|=3\cdot 2^{d-2}\) is such that \(V=V^{1,0}\cup V^{1,1}\) [\hyperlink{ref-CR9}{9}]. Then, by induction on \textit{d}, thanks to the structure of \textit{V} we can get a slightly better estimation: \(S_d<2^d+2^{d-1}\).

We also conjecture that if \(0\le N\le 2^{d-2}\) and \(V\) is a neighborly \(d\)-code without twin pairs with \(|V|=2^d-N\), then, up to isomorphism, \(V=V^{1,0}\cup V^{1,1}\) (we are able verify this conjecture for \(N\in \{0,1,2\}\)). Let us note that Corollary \hyperlink{FPar11}{1} is a small step in the direction of this conjecture: Corollary \hyperlink{FPar11}{1} says that if \textit{M} is rather small, then \(|V^{1,*}|\) is rather small, and \(|V^{1,0}|,|V^{1,1}|\) are rather big.

An inflation of a code \(V\subset A^n\) was a very useful tool in the proof of Theorem \hyperlink{FPar2}{2} and we believe that also other interesting properties of such codes can be discovered thanks to it.

\section*{Acknowledgements}

The author thanks Tony Newton for pointing out that Corollary~1 and the second conjecture in Remarks~2 omitted the assumption that the code has no twin pairs. This assumption was intended in both statements.

\bigskip\noindent Andrzej P. Kisielewicz\\
Wydzia{\l} Matematyki, Informatyki i Ekonometrii\\
Uniwersytet Zielonog\'orski\\
ul. Podg\'orna 50\\
65-246 Zielona G\'ora\\Poland\\
e-mail: \texttt{A.Kisielewicz@wmie.uz.zgora.pl}

{\small\noindent Revised version of the article published in
\href{https://doi.org/10.1007/s00026-024-00694-2}{\textit{Annals of Combinatorics} \textbf{28} (2024), 733--748}.
The published article is licensed under \href{https://creativecommons.org/licenses/by/4.0/}{CC BY 4.0}.\par}
\end{document}